\documentclass[11pt]{article}
\usepackage{graphicx}
\usepackage[all]{xy}
\usepackage{latexsym,amssymb,amsmath,epsfig}
\usepackage[usenames,dvipsnames]{pstricks}
\usepackage{epsfig}
\usepackage{epstopdf} 
\usepackage{epstopdf}
\usepackage{graphicx}
\usepackage{pstricks}
\usepackage{epsfig}
\usepackage{epstopdf}
\usepackage{epsfig}
\usepackage{amsthm}
\usepackage{cite}
\usepackage{caption}
\usepackage{authblk}
\usepackage{mathtools}
\usepackage[T1]{fontenc}
\usepackage[utf8]{inputenc}
\usepackage{tabularx,ragged2e,booktabs,caption}
\newcolumntype{C}[1]{>{\Centering}m{#1}}

\usepackage{lscape}
\usepackage{pdflscape}
\usepackage{float,lscape}
\usepackage{qtree}

\author[1]{Ali Raza}
\author[1]{Sonia Naseer}
\author[1]{FD Zaman}
\author[2]{AH Kara}
\affil[1]{Abdus Salaam School of Mathematical Sciences, Government College University, Lahore, Pakistan}
\affil[2]{School of Mathematics, University of the Witwatersrand, Wits 2050, South Africa}

\begin{document}

\title{Optimal System and Conservation Laws for the Generalized Fisher Equation in Cylindrical Coordinates}
\maketitle


\begin{abstract} The reaction diffusion equation arises in physical situations in problems from population growth, genetics and physical sciences. We consider the generalised Fisher equation in cylindrical coordinates from Lie theory stand point. An invariance method is performed and the optimal set of nonequivalent symmetries is obtained. Finally, the conservation laws are constructed using 'multiplier method'. We determine multipliers as functions of the dependent and independent variables only. The conservation laws are computed and presented in terms of conserved vector corresponding to each multiplier. \end{abstract}


\section{Introduction}
The Fisher equation for population dynamics was first proposed by Fisher in his article on advantageous genes in 1937 \cite{b1}.
\begin{equation}
u_t - u_{xx} = u(1-u),
\end{equation}
where, $ 0 \leq u(x,t) \leq 1$. Some of the history of the Fisher equation and its many forms are described by Rosa, Bruzón and Gandarias in reference \cite{b2}. The coupling of kinetics and diffusion gives rise to one dimensional Fisher equation whose scalar case can be written as
\begin{equation}
u_t - D u_{xx} = f(u),
\end{equation}
where $u$, $f(u)$ and $D$ represent the concentration, kinetics and diffusion coefficient respectively. Here the diffusion term $D$ is taken as constant. A simple case of nonlinear reaction diffusion equation  equation is \cite{b3}
\begin{equation}
u_t - D u_{xx} = k u (1-u).
\end{equation}
The Fisher equation belongs to the class of reaction diffusion equations that models population growth, wave propagation of advantageous genes, reacting diffusion process and other biological systems.\\
Generalised form of the Fisher equation and its modifications in the Cartesian coordinates as well as cylindrical case has attracted the attention of researchers. Most of the physical and biological models are described by nonlinear partial differential equations. In the one-dimensional case, the density dependent diffusion equation which is called generalised Fisher equation is investigated by Murray \cite{b3} which is given by
\begin{equation} \label{eq1}
u_t-\big(h(u)u_x)_x = f(u),
\end{equation}
where $h(u)$ is diffusion coefficient depending on $u(x,t)$. Rosa, Bruzón and Gandarias determined the Lie symmetries, travelling wave solutions, optimal systems and conservation laws using  Ibragimov \cite{b4} approach. In another article, Rosa and Gandarias \cite{b5} presented conservation laws using multiplier approach proposed by Anco and Bluman \cite{b6,b7}.  Bokhari, et. al. in \cite{b8} found Lie point symmetries of generalized Fisher equation (\ref{eq1}) and presented some invariant solutions for some cases. \\
In higher dimensions, more generalized form of Fisher equation can be written as
\begin{equation} \label{eq2}
u_t-\nabla \cdot(f(u)\nabla u)=g(u).
\end{equation}
To fit a particular model, we can write Fisher equation in cylindrical coordinates as many practical problems desire the use of cylindrical coordinates. By assuming radial symmetry the equation (\ref{eq2}) reduces to
\begin{equation} \label{eq3}
u_t-\frac{1}{x}  \cdot(xf(u)u_x )_x=g(u),
\end{equation}
which is the generalised Fisher equation in cylindrical coordinates denoted by $R[u]$. For particular case, by considering the values of $f(u)=u, g(u)=u(1-u)$, equation (\ref{eq3}) becomes
\begin{equation} \label{eq4}
u_t-\frac{1}{x} (xu  u_x )_x=u(1-u).
\end{equation}
Bokhari, et. al. performed a Lie symmetry analysis of equation (\ref{eq4}) and found some invariant solution in \cite{b9}.\\
O. O. Vaneeva, et. al. have considered the Fisher equation in which diffusivity and reaction term satisfy the power nonlinearities \cite{c1} and exponential nonlinearities \cite{c2} given by,
\begin{align}
f(x)u_t- (g(x) u^n u_x )_x &= h(x) u^m , \\
f(x)u_t- (g(x) e^{nu} u_x )_x &= h(x) e^{mu} .
\end{align}
The purpose of our study is to investigate generalised Fisher equation in cylindrical coordinates using Lie symmetry approach \cite{b10,b11} which is useful tool to find invariant exact solutions. The transformation that transform the given equation into simpler one is determined by the so called Lie symmetry generators that form Lie algebra under the commutator operation. The classification of Lie symmetry generators into non-equivalent classes is made easier with the help of the optimal system of subalgebras. The reduced set of ordinary differential equations and invariant solutions are obtained via reduction under the optimal system. Finally, certain conservation laws are discussed, which are useful in defining the geometrical and physical description of surfaces of partial differential equations generated by their solutions.

\section{Lie Point Symmetry Generators} \label{S1}
In this section we look at the group properties of equation (\ref{eq3}) including Lie symmetries and investigate several cases by considering function values of $f(u)$ and $g(u)$. The invariance generators in the space ($t,x,u$) admitted by equation (\ref{eq3}) are given by
\begin{equation}
X = \xi^1 \frac{\partial}{\partial t} + \xi^2 \frac{\partial}{\partial x} + \eta \frac{\partial}{\partial u},
\end{equation}
where $\xi^1$,$\xi^2$ and $\eta$ lie in $(t,x,u)-$space. The second order extended infinitesimal transformation $X^{[2]}$ and second order extended infinitesimals $\zeta_i$, $\zeta_{ij}$ that require to solve (\ref{eq3}) are defined in \cite{c3,b10}. The invariance criterion  for Lie symmetries is provided below, viz.,
\begin{equation} \label{E3}
X^{[2]} [u_t-\frac{1}{x}  \cdot(xf(u)u_x )_x=g(u) ]\Big|_{u_t-\frac{1}{x}  \cdot(xf(u)u_x )_x=g(u)}=0.
\end{equation}
The set of determining equations is obtained by expanding the equation (\ref{E3}), we have obtained
\begin{align}
\xi^2_x=0 ,  \ \ \xi^2_u=0, \ \  \xi^1_u &= 0, \\
 f_u[ -2 \xi^1_x + f_u \xi^2_t + \eta_u  ] + \eta f_{uu} + f \eta_{uu}  &=0, \\
f[ \xi^1 +x \xi^1_x - x \xi^2_t +x^2 \xi^1_{xx} - 2 x^2 \eta_{xu} ] \nonumber \\
- f_u [  x \eta + 2x^2 \eta_x ]  - x^2 \xi^1_{t} &= 0,   \\
  x \eta  g_u  +x g \xi^2_t - x  \eta_{t} - x g  \eta_u + f \eta_x + x f \eta_{xx}  &=0,  \\
   \eta f_u  - 2 f \xi^1_x + f  \xi^2_{t}  &=0.
\end{align}
Simplification lead us to the values of infinitesimals from the above set of determining equations $\xi^1=0$, $\xi^2=c$ and $\eta=0$. Thus for the arbitrary $f(u)$ and $g(u)$, we have
\begin{equation*}
	\xi^1(x,t,u)=0 \ \  , \ \  \xi^2(x,t,u)=c  \ \ , \ \ \eta(x,t,u)=0
\end{equation*}
Setting the value $c=1$ yields the symmetry generator for the equation (\ref{eq3}) for arbitrary values of $f(u)$ and $g(u)$. The symmetry generator we obtain forms the principal Lie algebra given by
\begin{equation*}
X_1 =   \frac{\partial}{\partial t}.
\end{equation*}
The following special cases of $f(u)$ and $g(u)$ are of interest:

\noindent\textbf{Case-1}. $f(u) = u$ and $g(u) = u(1-u)$ :\\
We have an additonal symmetry generator due to
 \begin{equation*}
X_2 =   e^{-t} \frac{\partial}{\partial t} + u e^{-t} \frac{\partial}{\partial u} .
\end{equation*}

\noindent\textbf{Case-2}. $f(u) = mu$ and $g(u) = pu^2$ :\\
Here, the symmetry algebra is extended by
 \begin{equation*}
X_2 =   t \frac{\partial}{\partial t} - u  \frac{\partial}{\partial u} .
\end{equation*}

\noindent\textbf{Case-3}. $f(u) = mu^n$ and $g(u) = pu^q$ :\\
A 2-dimensional algebra is provided with
 \begin{equation*}
X_2 =   (2qt-2t) \frac{\partial}{\partial t} + [qx-(n+1)x] \frac{\partial}{\partial x} -2u  \frac{\partial}{\partial u} .
\end{equation*}

\noindent\textbf{Case-4}. $f(u) = au^n$ and $g(u) = u$ :\\
Via $X_2$ and $X_3$, below, a three-dimensional algebra is given by the addtional symmetries,
\begin{align*}
X_2 &=   e^{-nt} \frac{\partial}{\partial t} + u e^{-nt}  \frac{\partial}{\partial u} ,  \\
X_3 &=   x \frac{\partial}{\partial x} + \frac{2}{n} u  \frac{\partial}{\partial u} .
\end{align*}

\noindent\textbf{Case-5}. $f(u) = \sqrt[n]{u} $ and $g(u) = u$ :\\
The algebra extended to three dimensional algebra by the addition of $X_2$ and $X_3$ given by,
\begin{align*}
X_2 &=   e^{-\frac{1}{n} t} \frac{\partial}{\partial t} + u e^{-\frac{1}{n} t}  \frac{\partial}{\partial u} ,  \\
X_3 &=  x \frac{\partial}{\partial x} + 2n u  \frac{\partial}{\partial u} .
\end{align*}

\noindent\textbf{Case-6}. $f(u) = au^n $ and $g(u) = \sqrt[m]{u}$ :\\
A two-dimensional algebra results with
\begin{align*}
X_2 &=   (2mt-2t) \frac{\partial}{\partial t}  + x[(n+1)m-1]\frac{\partial}{\partial x}+ 2mu  \frac{\partial}{\partial u}.
\end{align*}
In this study, we discuss the group properties by looking at different values of $f(u)$ and $g(u)$. It further will help us in constructing the one dimensional optimal system of subalgebras.

\section{One-dimensional Optimal System of Subalgebras} \label{S2}
The goal to find the optimal set of Lie symmetries is to obtain non-similar classes known as optimal system of subalgebras. Each element from optimal set represents the general class of symmetry and help in constructing the general class of invariant solution. To get this done, we follow the approach given by Olver \cite{b11} and demonstrated in \cite{c3}. The adjoint representation is given by
\begin{equation} \label{E5}
 Ad(\exp(\epsilon Y_i))(Y_j) = Y_j - \epsilon [{Y_i}, Y_j] + \frac{\epsilon^2}{2!} [{Y_i}, [{Y_i}, Y_j]] - \cdots,
\end{equation}
where $\epsilon \in \mathbb{R}$ and $[Y_i,Y_j]$ represents the Lie product defined by
\begin{equation}
[Y_i, Y_j] = Y_iY_j - Y_j Y_i.
\end{equation}

\noindent\textbf{The optimal system:}\\
The optimal system for all the cases discussed here:\\
\noindent\textbf{1. principal case}\\
For the principal algebra the optimal system of subalgebra is given by,
\begin{equation}
X^1=X_1.
\end{equation}

\noindent \textbf{2. Case-1 $f(u) = u$ and $g(u) = u(1-u)$ :}\\
The algebra is two dimensional which is represented as follows,
 \begin{equation*}
X_1 =   e^{-t} \frac{\partial}{\partial t} + u e^{-t} \frac{\partial}{\partial u} ,
X_2 =  \frac{\partial}{\partial t},
\end{equation*}
with non-zero commutator given by,
\begin{equation}
[X_1,X_2]=X_1 .
\end{equation}
The adjoint representations for case-1 presented in Table-\ref{tab:title2}.\\

\begin{minipage}{\linewidth}
\centering
\begin{tabularx}{\linewidth}{@{} C{1in} C{.85in} *4X @{}}\toprule[1.5pt]
\bf Ad & \bf $X_1$ & \bf $X_2$ \\
\bottomrule[1.25pt]
  $X_1$ & $X_1$ & $X_2 - \epsilon X_1 $    \\
\hline
$X_2$ & $ e^{\epsilon} X_1$ & $X_2$  \\
\bottomrule[1.25pt]
\end {tabularx}\par
\captionof{table}{Adjoint Table} \label{tab:title2}
\bigskip
\end{minipage}
With the action of the adjoint representation on general element $X=a_1 X_1 + a_2 X_2 \in \mathcal{L}_2$, we obtain the one-dimensional optimal system of all the  subalgebras for the case-1 given by \\
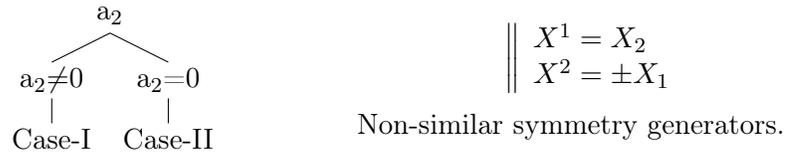
\begin{minipage}[t]{0.5\textwidth}
\Tree[ .a_2
		 [.a_2$\neq$0 {Case-I} ]
		 [.a_2=0   {Case-II}  ] ]
\end{minipage}
\begin{minipage}[t]{0.5\textwidth}
\[
   \left\|
                \begin{array}{ll}
                 X^1 =  X_2 \\
                  X^2 = \pm X_1  \\
                \end{array}
              \right.
  \]
Non-similar symmetry generators.
\end{minipage}\\
\\
The following are the complete details for each leaf:\\
\textbf{Case I.} $a_2\neq0$ with adjoint action on $X_1$, we have obtained
\begin{eqnarray}
X^{'}=Ad(e^{a X_1})X  =  a_2 X_2
\end{eqnarray}
\textbf{Case II.}  $a_2 =0$ with adjoint action on $X_2$, we possess
\begin{eqnarray}
X^{'}= Ad(e^{a X_2})X  = \pm X_1
\end{eqnarray}

\noindent \textbf{3. Case-2. $f(u) = mu$ and $g(u) = pu^2$}\\
The algebra along with their non-zero commutator for this case is given by,
 \begin{equation*}
X_1 =   \frac{\partial}{\partial t} ,
X_2 =  t \frac{\partial}{\partial t} - u \frac{\partial}{\partial u} \  ; \  \ [X_1,X_2]=X_1 .
\end{equation*}
For the given case, the one dimensional optimal system is same as in case-1 shown in Table-\ref{tab:title2}, since the algebra is same. \\
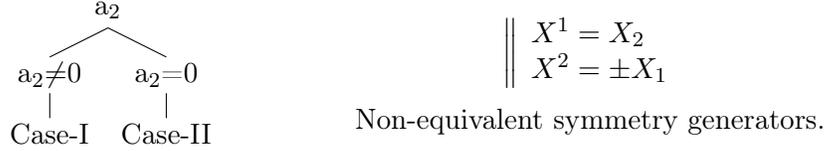
\begin{minipage}[t]{0.5\textwidth}
\Tree[ .a_2
		 [.a_2$\neq$0 {Case-I} ]
		 [.a_2=0   {Case-II}  ] ]
\end{minipage}
\begin{minipage}[t]{0.5\textwidth}
\[
   \left\|
                \begin{array}{ll}
                 X^1 =  X_2 \\
                  X^2 = \pm X_1  \\
                \end{array}
              \right.
  \]
Non-equivalent symmetry generators.
\end{minipage}\\
\\
\noindent \textbf{4. Case-3. $f(u) = mu^n$ and $g(u) = pu^q$ :}\\
The non-zero commutator for this case is given by,
\begin{equation}
[X_1,X_2]= 2(q-1)X_1 .
\end{equation}
The adjoint representations are shown in Table-\ref{tab:title3a},
when $f(u)=mu^n$ and $g(u)=pu^q $.\\

\begin{minipage}{\linewidth}
\centering
\begin{tabularx}{\linewidth}{@{} C{1in} C{.85in} *4X @{}}\toprule[1.5pt]
\bf Ad & \bf $X_1$ & \bf $X_2$ \\
\bottomrule[1.25pt]
  $X_1$ & $X_1$ & $X_2 - 2 \epsilon (q-1) X_1 $    \\
\hline
$X_2$ & $ e^{2 \epsilon (q-1)} X_1$ & $X_2$  \\
\bottomrule[1.25pt]
\end {tabularx}\par
\captionof{table}{Adjoint Table} \label{tab:title3a}
\bigskip
\end{minipage}
By the action of the adjoint representationgeneral on general element $X = a_1 X_1 + a_2 X_2 \in \mathcal{L}_2$ , we obtain the one-dimensional optimal system given by  \\
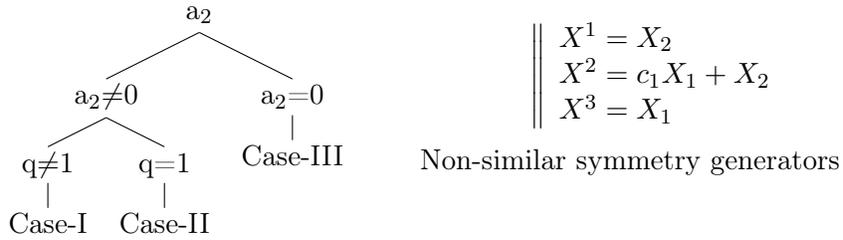
\begin{minipage}[t]{0.5\textwidth}
\Tree[ .a_2
		 [.a_2$\neq$0 [.q$\neq$1 {Case-I} ]
		 [.q=1   {Case-II}  ] ]
		 [.a_2=0   {Case-III}  ] ]
\end{minipage}
\begin{minipage}[t]{0.5\textwidth}
\[
   \left\|
                \begin{array}{ll}
                 X^1 =  X_2 \\
                  X^2 = c_1 X_1 + X_2  \\
                  X^3 =  X_1  \\
                \end{array}
              \right.
  \]
  Non-similar symmetry generators
\end{minipage}\\
\\
Here are given all of the details for each leaf:\\
\textbf{Case I.} $a_2\neq0$ and $q \neq 1$ . By the adjoint action on $X_1$, we possess
\begin{eqnarray}
X^{'}=Ad(e^{a X_1})X  =  a_2 X_2 = X_2
\end{eqnarray}
\textbf{Case II.} $a_2\neq0$ and $q=1$ . With the adjoint action on $X_1$ for $\epsilon = a$, we have
\begin{eqnarray}
X^{'}=Ad(e^{a X_1})X  =  a_1 X_1 + X_2
\end{eqnarray}
\textbf{Case-III.}  $a_2 =0$. For this case we apply adjoint action for any $Y$ and for any $\epsilon$, we obtain
\begin{eqnarray}
X^{'} = Ad(e^{ \epsilon Y})X  =  X_1
\end{eqnarray}

\noindent \textbf{5. Case-4. $f(u) = au^n$ and $g(u) = u$}\\
The algebra with basis $X_1, X_2$ and $X_3$ extend to three dimensional algebra as follows,
\begin{eqnarray*}
X_1 &=    \frac{\partial}{\partial t}, \\
X_2 &=  x \frac{\partial}{\partial x} +  \frac{2}{n} u \frac{\partial}{\partial u}  , \\
X_3 &=  e^{-nt} \frac{\partial}{\partial t} + u e^{-nt}   \frac{\partial}{\partial u}  ,
\end{eqnarray*}
with the non-zero commutator is given by,
\begin{equation}
[X_1,X_3]= -nX_3
\end{equation}
\\
The adjoint representations are given in Table-\ref{tab:title4a},
when $f(u)=mu^n$ and $g(u)=u $.\\

\begin{minipage}{\linewidth}
\centering
\begin{tabularx}{\linewidth}{@{} C{1in} C{.85in} *4X @{}}\toprule[1.5pt]
\bf Ad & \bf $X_1$ & \bf $X_2$ & \bf $X_3$ \\
\bottomrule[1.25pt]
  $X_1$ & $X_1$ & $X_2 $ & $ e^{n \epsilon} X_3  $    \\
\hline
$X_2$ & $  X_1$ & $X_2$ & $X_3$  \\
\hline
$X_3$ & $ X_1 - n \epsilon X_3 $ & $X_2$ & $X_3$  \\
\bottomrule[1.25pt]
\end {tabularx}\par
\captionof{table}{Adjoint Table} \label{tab:title4a}
\bigskip
\end{minipage}
The one-dimensional optimal system of all subalgebras for this case is determined by the adjoint representation, when applied on general element $X = a_1 X_1 + a_2 X_2 + a_3 X_3 \in \mathcal{L}_3$ given by, \\
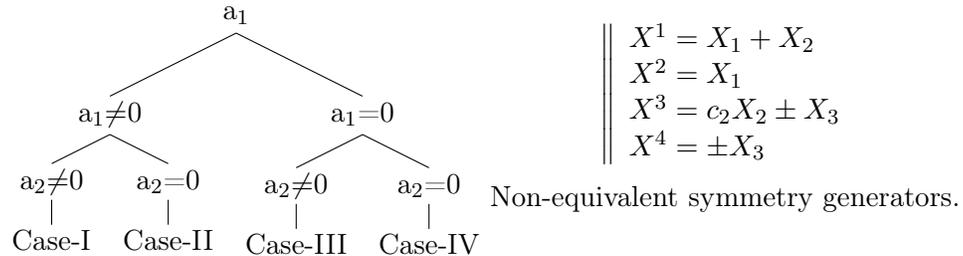
\begin{minipage}[t]{0.5\textwidth}
	\Tree[ .a_1
		 [.a_1$\neq$0 [.a_2$\neq$0 {Case-I} ]
		 [.a_2=0   {Case-II}  ] ]
		 [.a_1=0   [.a_2$\neq$0 {Case-III} ]
		 [.a_2=0   {Case-IV}  ]  ] ]
\end{minipage}
\begin{minipage}[t]{0.5\textwidth}

\[
   \left\|
                \begin{array}{ll}
                 X^1 =  X_1 + X_2 \\
                  X^2 =  X_1   \\
                  X^3 = c_2 X_2 \pm X_3  \\
                  X^4 =   \pm X_3
                \end{array}
              \right.
  \]
  Non-equivalent symmetry generators.
\end{minipage}\\
\\
Each leaf is described in detail as follows: \\
\textbf{Case I.} $a_1\neq0$ and $a_2 \neq 0$. We find the basis for optimal system using the adjoint action on $X_3$ for $\epsilon = a$,
\begin{eqnarray}
X^{'}=Ad(e^{a X_3})X  =  X_1 +  X_2.
\end{eqnarray}
\textbf{Case II.} $a_1\neq0$ and $a_2 = 0$. For this case, we obtain the basis for optimal system by the adjoint action on $X_3$ for $\epsilon = a$,
\begin{eqnarray}
X^{'}=Ad(e^{a X_3})X  =  X_1.
\end{eqnarray}
\textbf{Case III.}  $a_1=0$ and $a_2 \neq 0$. By the adjoint action on $X_1$ for $\epsilon = a$, we discover,
\begin{eqnarray}
X^{'}= Ad(e^{ \epsilon X_1 })X  = c_2 X_2 \pm X_3 ,
\end{eqnarray}
where $a=\frac{1}{n} ln|\pm \frac{1}{a_3}|$.\\
\textbf{Case IV.}  $a_1=0$ and $a_2 = 0$.  By the adjoint action on $X_1$ and for $\epsilon= a $, we own
\begin{eqnarray}
X^{'} = Ad(e^{ \epsilon X_1 })X  =  \pm X_3 .
\end{eqnarray}

\noindent \textbf{6. Case-5. $f(u) = \sqrt[n]{u} $ and $g(u) = u$}\\
In this case non-zero commutator is,
\begin{equation}
[X_1,X_2]= - \frac{1}{n} X_2 ,
\end{equation}
and the algebra is three dimensional where $X_1, X_2$ and $X_3$ are to be considered as follows,
\begin{eqnarray*}
X_1 &=    \frac{\partial}{\partial t}, \\
X_2 &=  e^{- \frac{1}{n} t} \frac{\partial}{\partial t} +  u e^{- \frac{1}{n} t} \frac{\partial}{\partial u}  , \\
X_3 &=  x \frac{\partial}{\partial x} + 2nu   \frac{\partial}{\partial u}  .
\end{eqnarray*}
The adjoint representations are described in Table \ref{tab:title5a},
when $f(u)= \sqrt[n]{u} $ and $g(u)=u $.\\

\begin{minipage}{\linewidth}
\centering
\begin{tabularx}{\linewidth}{@{} C{1in} C{.85in} *4X @{}}\toprule[1.5pt]
\bf Ad & \bf $X_1$ & \bf $X_2$ & \bf $X_3$ \\
\bottomrule[1.25pt]
  $X_1$ & $X_1$ & $ e^{\frac{1}{n} \epsilon} X_2 $ & $ X_3 $    \\
\hline
$X_2$ & $  X_1 - \frac{1}{n} \epsilon X_2 $ & $X_2$ & $X_3$  \\
\hline
$X_3$ & $   X_1$ & $X_2$ & $X_3$  \\
\bottomrule[1.25pt]
\end {tabularx}\par
\captionof{table}{Adjoint Table} \label{tab:title5a}
\bigskip
\end{minipage}
The one-dimensional optimal system of all subalgebras is determined by the adjoint representation, when operated on general element $X=a_1 X_1 + a_2 X_2 + a_3 X_3 \in \mathcal{L}_3$ given by  \\
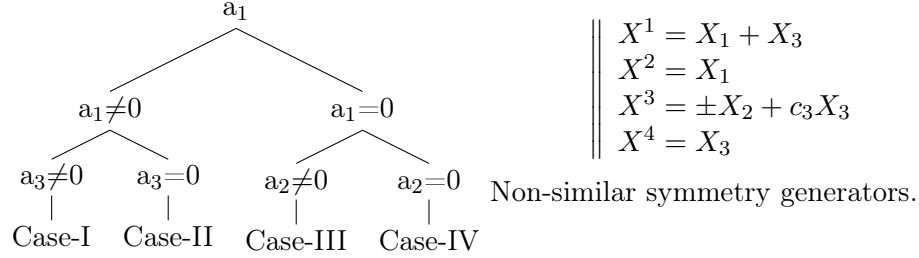
\begin{minipage}[t]{0.5\textwidth}
	\Tree[ .a_1
		 [.a_1$\neq$0 [.a_3$\neq$0 {Case-I} ]
		 [.a_3=0   {Case-II}  ] ]
		 [.a_1=0   [.a_2$\neq$0 {Case-III} ]
		 [.a_2=0   {Case-IV}  ]  ] ]
\end{minipage}
\begin{minipage}[t]{0.5\textwidth}
\[
   \left\|
                \begin{array}{ll}
                 X^1 =  X_1 + X_3 \\
                  X^2 =  X_1   \\
                  X^3 = \pm X_2 + c_3 X_3  \\
                  X^4 =    X_3
                \end{array}
              \right.
  \]
  Non-similar symmetry generators.
\end{minipage}\\
\\
Each leaf is discussed in detail as follows:\\
\textbf{Case I.} $a_1\neq0$ and $a_3 \neq 0$ . By the adjoint action on $X_2$, we possess
\begin{eqnarray}
X^{'}=Ad(e^{a X_2})X  =  X_1 +  X_3.
\end{eqnarray}
\textbf{Case II.} $a_1\neq0$ and $a_3= 0$. With the adjoint action on $X_2$,  we acquire
\begin{eqnarray}
X^{'}=Ad(e^{a X_2})X  =  X_1.
\end{eqnarray}
\textbf{Case III.}  $a_1=0$ and $a_2 \neq 0$. As a result of an adjoint action on $X_1$ and for $\epsilon= a $, we have received
\begin{eqnarray}
X^{'} = Ad(e^{ \epsilon X_1 })X  = \pm X_2 + c_3 X_3 ,
\end{eqnarray}
where $a=(n) ln|\pm \frac{1}{a_2}|$.\\
\textbf{Case IV.}  $a_1=0$ and $a_2 = 0$. The adjoint action on $X_1$ and for $\epsilon= a $ provide
\begin{eqnarray}
X^{'} = Ad(e^{ \epsilon X_1 })X  =   X_3 .
\end{eqnarray}

\noindent \textbf{Case-6. $f(u) = au^n $ and $g(u) = \sqrt[m]{u}$}\\
The algebra is two dimensional algebra and basis $X_1$, $X_2$ are represented as follows,
 \begin{equation*}
X_1 =    \frac{\partial}{\partial t} ,
X_2 =  2t(m-1) \frac{\partial}{\partial t} + x[(n+1)m-1] \frac{\partial}{\partial x} +  2mu \frac{\partial}{\partial u}  ,
\end{equation*}
with non-zero commutator,
\begin{equation}
[X_1,X_2]= 2(m-1)X_1 .
\end{equation}
The adjoint representations are listed in Table \ref{tab:title6a},
when $f(u)=au^n$ and $g(u)= \sqrt[m]{u} $.\\

\begin{minipage}{\linewidth}
\centering
\begin{tabularx}{\linewidth}{@{} C{1in} C{.85in} *4X @{}}\toprule[1.5pt]
\bf Ad & \bf $X_1$ & \bf $X_2$ \\
\bottomrule[1.25pt]
  $X_1$ & $X_1$ & $X_2 - 2 \epsilon (m-1) X_1 $    \\
\hline
$X_2$ & $ e^{2 \epsilon (m-1)} X_1$ & $X_2$  \\
\bottomrule[1.25pt]
\end {tabularx}\par
\captionof{table}{Adjoint Table} \label{tab:title6a}
\bigskip
\end{minipage}
With the action of the adjoint representation on general element $X =a_1 X_1 + a_2 X_2 \in \mathcal{L}_2$, the one-dimensional optimal system of all the  subalgebras for this case are given below  \\
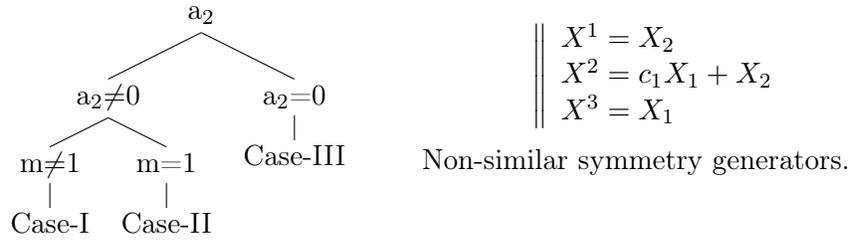
\begin{minipage}[t]{0.5\textwidth}
	\Tree[ .a_2
		 [.a_2$\neq$0 [.m$\neq$1 {Case-I} ]
		 [.m=1   {Case-II}  ] ]
		 [.a_2=0   {Case-III}  ] ]
\end{minipage}
\begin{minipage}[t]{0.5\textwidth}
\[
   \left\|
                \begin{array}{ll}
                 X^1 =  X_2 \\
                  X^2 = c_1 X_1 + X_2  \\
                  X^3 =  X_1  \\
                \end{array}
              \right.
  \]
  Non-similar symmetry generators.
\end{minipage}\\
\\
Details for each leaf are provided below:\\
\textbf{Case I.} $a_2\neq0$ and $m \neq 1$ . We accomplish this by the adjoint action on $X_1$ for $\epsilon = a$,
\begin{eqnarray}
X^{'}=Ad(e^{a X_1})X  =  a_2 X_2 = X_2
\end{eqnarray}
\textbf{Case II.} $a_2\neq0$ and $m=1$ . By the adjoint action on $X_1$ for $\epsilon = a$, we have obtained
\begin{eqnarray}
X^{'}=Ad(e^{a X_1})X  =  a_1 X_1 + X_2
\end{eqnarray}
\textbf{Case III.}  $a_2 =0$. The adjoint action for any $Y$ and for any $\epsilon$ results in
\begin{eqnarray}
X^{'} = Ad(e^{ \epsilon Y})X  =  X_1
\end{eqnarray}

\section{Reduction and Similarity Solutions} \label{S3}

The optimal system of subalgebras from section-\ref{S2} will be used to reduce the equation-(\ref{eq3}). The procedure is well-known and it is outlined in \cite{b10,b11,b12} and demonstrated in \cite{c3}. The characteristic equation is given by

\begin{equation}
\frac{dt}{\xi^1(t,x,u)} = \frac{dx}{\xi^2(t,x,u)} = \frac{du}{\eta(t,x,u)}.
\end{equation}
Complete reduction is presented for each case that was investigated in this study for different value of $f(u)$ and $g(u)$. Reductions of equation (\ref{eq3}) for a nontrivial symmetry generator from the optimal system (section \ref{S2}) is given as follows:\\
For instance, we consider $X^{2}$ from Case-3.
\begin{equation*}
X^{2} = c_1 X_1 + X_2  ,
\end{equation*}
and its respective invariant variables are appeared as,
\begin{align*}
u &= F(\alpha) x^{\frac{2}{n+1}}  \\
  \alpha &= \frac{1}{2} \frac{(2qt+c-2t)}{q-1} x^{\frac{2(q-1)}{n+1}}.
\end{align*}
thus the invariant solution is computed as,
\begin{equation} \label{E7}
u = F(\frac{1}{2} \frac{(2qt+c-2t)}{q-1} x^{\frac{2(q-1)}{n+1}}) x^{\frac{2}{n+1}} .
\end{equation}
Therefore, (\ref{E7}) represents the solution of the equation (\ref{eq3}) if $F$  satisfies the ODE,
\begin{eqnarray*}
 4m F^{n-1} [ \alpha^2 F (q-1)^2 F^{''} + n \alpha^2 (q-1)^2 {F^{'}}^2  ,\\
      + (q+2n+1) \alpha (q-1) F F^{'}     + F^2 (n+1) ]=0  .
\end{eqnarray*}
The complete set of reduced ODEs under optimal system are presented in Table \ref{tabler1},\ref{tabler2},\ref{tabler3},\ref{tabler4}, and \ref{tabler5}.

\begin{center}
\begin{tabular}{ l | l | l }\toprule[1.5pt]
 \bf  Generators & \bf Reduced ODE & \bf Similarity Solutions  \\
\midrule[1.5pt]
\bf & For Principal case. \\
 \midrule[1.5pt]
   $X^{1}=X_1 $ & $F^{''} - \frac{1}{\alpha} F^{'} + \frac{1}{\alpha} F(\alpha) = 0  $ &  $ u = F(\alpha)$ \\
   & & $ \alpha = x$    \\
\midrule[1.5pt]
\end{tabular}
\captionof{table}{Reduction set of ODEs} \label{tabler1}
\end{center}

\begin{center}
\begin{tabular}{ l | l | l }\toprule[1.5pt]
 \bf  Generators & \bf Reduced ODEs & \bf Similarity Solutions  \\
\midrule[1.5pt]
\bf For Case 1. \\
 \midrule[1.5pt]
    $X^{1} = X_2 $ & $ \alpha F^2 - \alpha F F^{''} - \alpha {F^{'}}^2 - F F^{'}=0$    &  $ u = F(\alpha) e^t $\\ & & $ \alpha = x$ \\
   \hline
    $X^{2}= X_1 $ &  $ F  F^{'} + \alpha  {F^{'}}^2  + \alpha F  F^{''} + \alpha F - \alpha F^2  =0 $
    & $ u = F(x)$ \\
    & &  $ \alpha = x$ \\
    \midrule[1.5pt]
\bf  For Case 2. \\
 \midrule[1.5pt]
    $X^{1} = X_1 $ &  $  m {F^{'}}^2 + m F F^{''} + \alpha p F^2 = 0   $
    & $ u = F(x)  $ \\
    & &  $ \alpha = x$ \\
   \hline
    $X^{2}= X_2 $ &  $ \alpha p  F^2 + m F   F^{''} + m {F^{'}}^2 - \alpha F   =0 $
   & $ u = \frac{F(x)}{t}$ \\
   & &  $ \alpha = x$ \\
   \midrule[1.5pt]
\bf  For Case 3. \\
 \midrule[1.5pt]
    $X^{1} = X_2 $ &  $  4 \alpha^{\frac{n}{2(q-1)}} m F^{n-1} $ \\
    &  $ [ (q+2n+1)(q-1) F F^{'} \alpha^{\frac{q}{q-1}} \alpha^{-\frac{1}{2(q-1)}} $ \\
    & $ + \alpha^{\frac{2q}{q-1}}(q-1)^2 [{F^{'}}^2 n + F^{''}F] $ \\
    & $\alpha^{\frac{-3}{2(q-1)}}+ \alpha^{\frac{1}{2(q-1)}}(n+1) F^2 ] = 0   $
    &  $ u = F(t x^{2\frac{q-1}{n+1}}) x^{\frac{2}{n+1}}  $ \\
    & & $ \alpha = t x^{2\frac{q-1}{n+1}}$             \\
   \hline
    $X^{2}= c X_1 + X_2 $ &  $ 4m F^{n-1} [ \alpha^2 F (q-1)^2 F^{''} + n \alpha^2 (q-1)^2 {F^{'}}^2 $\\
    & $ + (q+2n+1) \alpha (q-1) F F^{'} + F^2 (n+1) ]=0  $
    &   $ u = F(x) x^{\frac{2}{n+1}}$  \\
    & & $ \alpha = \frac{1}{2} \frac{(2qt+c-2t)}{q-1} x^{\frac{2(q-1)}{n+1}}$            \\
    \hline
    $X^{3}= X_1  $ &  $ mF^n F^{'} + \alpha m n F^{n-1}  {F^{'}}^2 - \alpha m F^n F^{''} - \alpha p F^{q} =0  $
    & $ u = F(x) $ \\
    & &  $ \alpha = x$ \\
\midrule[1.5pt]
\end{tabular}
\captionof{table}{Reduction set of ODEs} \label{tabler2}
\end{center}

\begin{center}
\begin{tabular}{ l | l | l }\toprule[1.5pt]
 \bf  Generators & \bf Reduced ODEs & \bf Similarity Solutions \\
\midrule[1.5pt]
\bf For Case 4. \\
 \midrule[1.5pt]
  $X^{1} = X_1 + X_2 $ &  $  n^2 F^n  F^{''} - F^{'}[4n^2 F^n + 4nF^n +n^2] $ \\
    &  $ +[ n^2  +4nF^n]F + 4 F^{n+1} - n^3 {F^{'}}^2 F^{n-1} = 0 $
    &   $ u = F(- ln(x) + t) x^{2/n}  $ \\
    & &     $ \alpha = - ln(x) + t $ \\
   \hline
    $X^{2}= X_1 $ &  $ F^n F^{'} - \alpha n F^{n-1} {F^{'}}^2 -  \alpha  F^n F^{''} - \alpha F =0  $
    &  $ u = F(x) x^{\frac{2}{n+1}}$ \\
    & &  $ \alpha = \frac{1}{2} \frac{(2qt+c-2t)}{q-1} x^{\frac{2(q-1)}{n+1}}$ \\
    \hline
    $X^{3}= X_2 + X_3   $ &  $ F^n F^{''} + \frac{1}{n^2} [ n^3 {F^{'}}^2 -4n(n+1) F F^{'} $ \\
    & $+ 4(n+1) F^2 ] + F^{n-1} -n  F^{'} =0  $
    & $ u = F(\frac{-n ln(x)- e^{nt}}{n}) n^{\frac{-1}{n}} e^t $ \\
    & &  $ \alpha = \frac{e^{nt} -n ln(x) }{n}$ \\
    \hline
    $X^{4}=  X_3   $ &  $ [\alpha F^n] F^{''} +  F^n F^{'} + \alpha n F^{n-1} {F^{'}}^2 =0  $
    &  $ u = F(\alpha) e^{-t} $ \\
    & &  $ \alpha = x$ \\
\midrule[1.5pt]
\end{tabular}
\captionof{table}{Reduction set of ODEs} \label{tabler3}
\end{center}

\begin{center}
\begin{tabular}{ l | l | l }\toprule[1.5pt]
 \bf  Generators & \bf Reduced ODEs & \bf Similarity Solutions \\
\midrule[1.5pt]
\bf For Case 5. \\
 \midrule[1.5pt]
    $X^{1} = X_1 + X_3 $ &  $  [nF^{''} +   {F^{'}}^2 F^{-1} - 4n (1+n) F^{'} + $ \\
    &  $ 4n^2 (1+n) F ]  F^{1/n} + n[ F - F^{'}  ] =0 $
    &   $ u = F(-ln(x) + t) x^{2n}  $  \\
    & & $ \alpha = -ln(x) + t $ \\
     \hline
    $X^{2}=  X_1 $ &  $ n F^{\frac{1}{n}} F^{'} - \alpha F^{\frac{(1-n)}{n}}   {F^{'}}^2 $ \\
    & $ - \alpha n F^{1/n} F^{''} - \alpha n F =0  $
   & $ u = F(\alpha) $ \\
   &  &   $ \alpha = x$ \\
     \hline
    $X^{3}= X_2 + c X_3   $ , $c=1$ &  $ (n^n e^{-\alpha (1-2n)}) [ n F^{''} + {F^{'}}^2 F^{-1}  $ \\
    & $- 4n (1+n) F^{'} $\\
    & $+ 4n^2 F (1+n ) ] F^{1/n} - F^{'} ] = 0  $
    &  $ u = F(e^{\frac{t}{n}} n - ln(x)) x^{2n} e^t n^n  $ \\
    & & $ \alpha = e^{\frac{t}{n}} n - ln(x) $ \\
    \hline
    $X^{4}=  X_3   $ &  $ 4 n^2 F F^{1/n}  + 4n F F^{1/n}  - F^{'} + F  = 0  $
    &  $ u = F(\alpha) x^{2n}   $ \\
    & &  $ \alpha = t $ \\
\midrule[1.5pt]
\end{tabular}
\captionof{table}{Reduction set of ODEs} \label{tabler4}
\end{center}

\begin{center}
\begin{tabular}{ l | l | l }\toprule[1.5pt]
  \bf  Generators & \bf Reduced ODEs & \bf Similarity Solutions \\
\midrule[1.5pt]
\bf For Case 6. \\
   \midrule[1.5pt]
    $X^{1}=  X_2   $ &  $  \alpha^{\frac{-mn}{2m-2}} [  [ ((n+1)m - 1)^2[   {F^{'}} - F^{1/m} ] $\\
     &$  - 4m (n+1) F^n F^2   ]\cdot  \alpha^{\frac{-1}{2m-2}} \alpha^{\frac{-m}{ 2m-2 }} $\\
    &$ + [\alpha^{\frac{-5}{2m-2}} \alpha^{\frac{3m}{2m-2}} (m-1)( n{F^{'}}^2 + F {F^{''}}  ) F^{-1} $ \\
    & $ - 2 {F^{'}} \alpha^{\frac{m}{2m-2}} \alpha^{\frac{-3}{2m-2}} (\frac{1}{2} + (n+\frac{1}{2})m )  ](m-1) F^n ] =0   $ & $ u = F(\alpha) x^{\frac{2m}{mn+m-1}}    $       \\
    &  &  $ \alpha = t x^{\frac{-2(m-1)}{mn+m-1}} $ \\
       \hline
    $ X^{2}= c_1 X_1 +  X_2  $ &  $ [(n+1)m - 1]^2 F^{1/m} + 4 F^n \alpha^2 (m-1)^2 F^{''} $ \\
    & $ + 4 F^{n-1} n \alpha^2 (m-1)^2 {F^{'}}^2 $    \\
    &    $ -8 \alpha (m-1) [\frac{1}{2} +(n+\frac{1}{2})m   ] F^n  {F^{'}} $ \\
    & $- [(n+1)m-1 ]^2  {F^{'}} + F^n F m^2 (n+1)  =0  $ & $ u = F(\alpha)x^{\frac{2m}{mn+m-1}} $ \\
    &  & $ \alpha = \frac{1}{2}  \frac{(2mt-2t+1) x^{\frac{-2(m-1)}{mn+m-1}})}{m-1}  $ \\
    \hline
    $X^{3}=  X_1   $ &  $ F^5 {F^{'}} + 5 \alpha F^4 {F^{'}}^2 + \alpha F^5 {F^{''}} + \alpha F^{\frac{1}{m}} =0   $  & $ u = F(\alpha)    $ \\
& &  $ \alpha = t $ \\
\midrule[1.5pt]
\end{tabular}
\captionof{table}{Reduction set of ODEs} \label{tabler5}
\end{center}

The complete set of ODEs corresponding to each case under optimal system is presented in this section. In the next session we will look for some conservation laws.

\section{Conservation Laws}
The study of conservation laws for the partial differential equations is significant because they explain the geometrical properties of the surface formed by solution and are simply mathematical objects. Conservation laws via multiplier approach are presented by Rosa and Gandarias in \cite{b2}. Rosa and Gandarias in \cite{b2} investigated the conservation laws of the generalised Fisher equation (refeq1) in three cases, each with distinct values of $f(u)$ and $h(u)$ and computed some multipliers and corresponding conservation laws with no dependence on variable $u_t$.\\
We follow the direct approach \cite{a1,a2,a3,a4}  to find conservation laws presented in Anco and Bluman \cite{b6,b7}. In this work, we use the multiplier approach to compute the conservation laws of the generalised Fisher equation in cylindrical coordinates for different value of $f(u)$ and $g(u)$. Multipliers are dependent on variables \{$x,t,u$\} for which the corresponding conservation laws are formulated. Consider a multiplier $\Lambda(x,t,u)$  satisfying
\begin{equation}\label{euler}
E_{u}(\Lambda R[u])=0,
\end{equation}
where $E_{u}$ represents the Euler operator given as
\begin{equation}
E_{u}=\frac{\partial }{\partial u}-D_t\frac{\partial}{\partial u_t}+D_xD_x\frac{\partial}{\partial u_{xx}}+D_xD_t\frac{\partial}{\partial u_{xt}}+D_tD_t\frac{\partial}{\partial u_{tt}}\cdots,
\end{equation}
and $D_x, D_t$ are total derivative operators and can be represented as,
\begin{equation*}
D_i=\frac{\partial}{\partial x^i}+u_i^\alpha\frac{\partial}{\partial u^\alpha}+u_{ij}^\alpha\frac{\partial}{\partial u_j^\alpha}+\cdot\cdot\cdot~~i=1,2,\cdots.
\end{equation*}
A conserved vector is a conserved quantity of (\ref{eq3}) which is a $2$-tuple $T=(T^t,T^x)$ such that it holds for all solutions of equation-(\ref{eq3}) given as,
\begin{equation} \label{11}
D_t T^t + D_x T^x = 0 .
\end{equation}
The multipliers $ \Lambda \{ x,t,u \} $ of the PDE (\ref{eq3}) has the property
\begin{equation} \label{E12}
D_tT^t(x,t,u)+D_xT^xt(x,t,u)=\Lambda(x,t,u)R[u] .
\end{equation}
where multipliers $\Lambda$ are depending independent variables on $x, \ t$ and dependent variable $  u $.s

\subsection{Multipliers and Corresponding Conservation Laws}

Consider the multipliers by $\Lambda(x,t,u)$ and expanding the expression (\ref{E12}) which is the invariance criteria to find conservation laws, leads us to the expression helps in finding the set of determining equations as follows:
\begin{eqnarray}
-xf'(u)\Lambda_uu^2_x-xf(u)\Lambda_uu_{xx}-xg(u)\Lambda_u-xg'(u)\Lambda \nonumber  \\
+f(u)\Lambda_x+2xf'(u)\Lambda_uu_x^2-x\Lambda_t-2f(u)\Lambda_x-xf(u)\Lambda_{xx}\\
-2f(u)\Lambda_uu_x-xf(u)\Lambda_{uu}u_x^2-2xf(u)\Lambda_{ux}u_x^2-xf(u)\Lambda_uu_{xx}=0. \nonumber
\end{eqnarray}
The above equation is an algebraic equation in terms of $u$, its derivatives, products and powers of derivatives. We compare like terms to earn system of following determining equations
\begin{eqnarray}
xf(u)\Lambda_u=0
xf(u)\Lambda_{uu}+xf'(u)\Lambda_u &= 0  \label{De1}   \\
 2xf(u)\Lambda_{ux}+f(u)\Lambda_u &= 0  \label{De2}  \\
xf(u)\Lambda_{xx}+f(u)\Lambda_x+xg(u)\Lambda_u+x\Lambda_t+xg'(u)\Lambda &= 0 \label{De3} .
\end{eqnarray}
By solving above system, for different values of $f(u)$ and $g(u)$, we obtain different multipliers as follows:\\

\noindent\textbf{Conservation laws}\\
Conservation laws for each case are discussed here:\\
\noindent\textbf{1. Case-1.}\\
The set of equations (\ref{De1},\ref{De2},\ref{De3}) that follow take the form:
\begin{eqnarray}
-x\Lambda_uu^2_x-xu\Lambda_uu_{xx}-x\big(u(1-u)\big)\Lambda_u-x(1-2u)\Lambda \nonumber
+u\Lambda_x+2x\Lambda_uu_x^2-x\Lambda_t \\
-2u\Lambda_x-xu\Lambda_{xx} \nonumber -2u\Lambda_uu_x-xu\Lambda_{uu}u_x^2-2xu\Lambda_{ux}u_x^2-xu\Lambda_uu_{xx}=0. \nonumber
\end{eqnarray}
After comparing like terms and solving system of determining equations, we obtain the multiplier
\begin{equation} \label{case1}
\Lambda=c_1e^{-t}\big[I_0(\sqrt{2}x)\big]+c_2e^{-t}\big[K_0(\sqrt{2}x)\big],
\end{equation}
where $I_0(\sqrt{2}x)$ and $K_0(\sqrt{2}x)$ denote modified Bessels functions of first and second kind of order zero, respectively.\\
From the equation (\ref{case1}) for multiplier $\Lambda_1=e^{-t}\big[I_0(\sqrt{2}x)\big]$, equation (\ref{E12}) leads us to the characteristic equation corresponding to multiplier $\Lambda_1$:
\begin{eqnarray*}
T^t_t + T^t_u u_t +T^t_{u_{x}} u_{tx} + T^x_x + T^x_uu_x + T^x_{u_{x}} u_{xx}  \\
= e^{-t} [ I_0 (\sqrt[2]{x}) ](xu_t - uu_x - xu^2_x - xuu_{xx} -xu+xu^2)
\end{eqnarray*}
After comparing the terms of highest order $u_{tx},u_{xx}$, we obtain the set of determining equations:
\begin{eqnarray}
T^x_{u_x}=-xue^{-t}I_0(\sqrt{2}x) \\
T^t_{u_x}=0
\end{eqnarray}
By solving the above set of equations, we obtain the set of local conservation laws of (\ref{eq3}) for the case when $f(u)=u$ and $ g(u)=u(1-u)$ which are presented in the form of conserved vectors $T_n = (T^t , \  \ T^x )$, where $n=1,2$:
\begin{eqnarray}
T_1 = ( xue^{-t}[I_0(\sqrt{2}x)] ,  \  \ \frac{1}{2}(\sqrt{2}xu^2e^{-t}[I_1(\sqrt{2}x)]-xue^{-t}[I_0(\sqrt{2}x)]u_x)),
\end{eqnarray}
\begin{eqnarray}
T_2 =
( xue^{-t}[K_0(\sqrt{2}x)] , \  \ -\frac{1}{2}(\sqrt{2}xu^2e^{-t}[K_1(\sqrt{2}x)]-xue^{-t}\big[K_0(\sqrt{2}x)\big]u_x) ).
\end{eqnarray}
which corresponds to the multipliers $\{\Lambda_1,\Lambda_2\} $ which are given by $\{e^{-t}[I_0(\sqrt{2}x)],$ $ \ \ e^{-t}[K_0(\sqrt{2}x)]\}$ respectively. Following the same procedure we find the conserved vectors corresponding to each multipliers for the rest of the cases as well.  \\

\noindent\textbf{2. Case-2}\\
The multipliers for the case when we consider the following values of function $f(u) = au $ and $g(u)= pu^2$ are given by:
\begin{equation}
\Lambda=c_1[J_0(\sqrt{2}x)]+c_2[Y_0(\sqrt{2}x)].
\end{equation}
For the multiplier $\Lambda_1=[J_0(\sqrt{2}x)]$, we deduce the conserved vectors:
\begin{align}
T_1 = (T^t , \  T^x) =   ( xu[J_0(\sqrt{2}x)], -\frac{1}{2}(\sqrt{2}xu^2\big[J_1(\sqrt{2}x)\big])-xu\big[J_0(\sqrt{2}x)\big]u_x   ).
\end{align}
For the multiplier $\Lambda_2=\big[Y_0(\sqrt{2}x)\big]$, we find
\begin{align}
T_2 = (T^t , \  T^x) =   ( xu[Y_0(\sqrt{2}x)],  -\frac{1}{2}(\sqrt{2}xu^2\big[Y_1(\sqrt{2}x)\big])-xu\big[Y_0(\sqrt{2}x)\big]u_x   ).
\end{align}

\noindent\textbf{3. Case-3}\\
The multipliers for the case when $f(u) = mu^n $ and $g(u)= pu^q$ are given by:
\begin{equation} \label{case5}
\Lambda=e^{-pqu^{q-1}t}\big[c_1+c_2ln(x)\big].
\end{equation}
For the multiplier $\Lambda_1=e^{-pqu^{q-1}t}$, we obtain the conservation laws given in the vector forms:
\begin{align}
T_1 = (T^t , \  T^x) =   ( xe^{-pqu^{q-1}t}, \ \ -xmu^ne^{-pqu^{q-1}t}u_x  ).
\end{align}\\
For the multiplier $\Lambda_2=e^{-pqu^{q-1}t}(lnx)$, we derive $T_2 = (T^t , \  T^x) $ which is given by
\begin{eqnarray}
T_2 =( xue^{-pqu^{q-1}t}(lnx), \ \ -xmu^ne^{-pqu^{q-1}t}(lnx)u_x+e^{-pqu^{q-1}t}\frac{mu^{n+1}}{m}  ).
\end{eqnarray}

\noindent\textbf{4. Case-4}
The multipliers for the case when $f(u)=mu^n, g(u)=u$ for all $n$, are given by,
\begin{equation} \label{case11}
\Lambda=e^{-t}\big[c_1+c_2ln(x)\big].
\end{equation}
Corresponding to multiplier $\Lambda_1=e^{-t}$, we acquire
\begin{align}
T_1 = (T^t , \  T^x) =   ( xue^{-t}, \ \ -xe^ue^{-t}u_x ).
\end{align}\\
Conservation law corresponding to the multiplier $\Lambda_1=e^{-t}(lnx)$, we possess
\begin{align}
T_2 = (T^t , \  T^x) =   ( xue^{-t}(lnx), \ \ -xe^ue^{-t}(lnx)u_x+e^x ).
\end{align}

\noindent\textbf{5. Case-5}
The multipliers for the case when $f(u)=\sqrt[n]{u}, g(u)=u$ for all $n$, are given by\\
\begin{equation} \label{case11}
\Lambda=e^{-t}[c_1+c_2ln(x)].
\end{equation}
For the multiplier $\Lambda_1=e^{-t}$, we obtain
\begin{align}
T_1 = (T^t , \  T^x) =   ( xue^{-t}, \ \ -x u_x e^{-t}  \sqrt[n]{u} ).
\end{align}\\
Conservation law corresponding to the multiplier $\Lambda_1=e^{-t}ln(x)$, given by
\begin{align}
T_2 = (T^t , \  T^x) =   ( xue^{-t}ln(x), \ \ - \sqrt[n]{u} e^{-t}( \frac{n}{n+1}u -  ln(x)x u_x ).
\end{align}\\
In a recent work, Anco \cite{b13}  proposed a study that establishes a connection between multipliers and the variational symmetry that underlies Noether's theorem. As previously stated, all of the conservation laws obtained by the multiplier approach have a clear application in mathematical analysis. Although some conservation laws describe physical quantities and explain natural phenomena like as momentum, energy, and object motion, the rest of the conservation laws describe the geometry of the surface formed by the partial differential equation.

\section{Conclusion} \label{S5}
The Lie symmetry analysis of the generalised Fisher equation in cylindrical coordinates (\ref{eq3}) was performed for arbitrary values of functions $f(u)$ and $g(u)$. One dimensional optimal system of subalgebras were classified and their detailed case-by-case study was presented through tree-leaf diagrams. Reduction has been performed and complete set of reduced ordinary differential equation under optimal system was listed in the Table \ref{tabler1}, \ref{tabler2}, \ref{tabler3}, \ref{tabler4} and Table \ref{tabler5}. Furthermore, we presented some conservation laws in the conserved form for the generalised Fisher equation in cylindrical coordinates.
\\
\\

\textbf{Ethics statement}: For this submitted article in  Acta Applicandae Mathematicae  titled `Optimal System and Conservation Laws for the Generalized Fisher Equation in Cylindrical Coordinates’ by Reza, Naseer, Zaman and Kara, we, the authors comply with the ethics of the journal.\\
No funding was received, the article is not submitted to another journal, all authors have made contributions to the article, for example.
\\
\\

\textbf{Conflict of Interest}: The authors declare that they have no conflict of interest.

 \end{document}